\documentclass[english]{article}
\usepackage[T1]{fontenc}
\usepackage[latin9]{inputenc}
\usepackage{babel}

\usepackage{array}
\usepackage{amsmath}
\usepackage{graphicx}
\usepackage{setspace}
\usepackage{amssymb}
\usepackage{esint}
\usepackage[unicode=true, pdfusetitle,
 bookmarks=true,bookmarksnumbered=false,bookmarksopen=false,
 breaklinks=false,pdfborder={0 0 1},backref=false,colorlinks=false]
 {hyperref}

\makeatletter

\providecommand{\tabularnewline}{\\}

\makeatother

\begin{document}

\title{\textbf{Good variation theory: a tauberian approach to the Riemann
Hypothesis }{\normalsize (Dedicated to the memory of Jonathan M. Borwein)}}

\author{Benoit Cloitre}
\maketitle
\begin{abstract}
In this note, I present a tauberian conjecture that I consider to
be the simplest and the best Tauberian reformulation of the Riemann
Hypothesis ($RH$) using good variation theory. 
\end{abstract}

\subsubsection*{Notations and definitions}

Let $\left(a(n)\right)_{n\geq1}$ be a sequence and $g$ be a function. 

\begin{flushleft}
We use the notations $A(n)=\sum_{k=1}^{n}a(k)$ and $A_{g}(n)=\sum_{k=1}^{n}a(k)g\left(\frac{k}{n}\right)$.
\par\end{flushleft}

\begin{flushleft}
The little Mellin transform of $g$ which is Riemann integrable on
$]0,1]$ is the analytic continuation of the function $g^{\star}$
defined for $\Re z<0$ by
\par\end{flushleft}

\[
g^{\star}(z):=\int_{0}^{1}g(t)t^{-z-1}dt\]

\begin{flushleft}
Next the analytic index of $g$ is defined by $\eta(g):=\min\left\{ \Re(\rho)\ \mid\ g^{\star}(\rho)=0\right\} $.
\par\end{flushleft}

\section*{Introduction}

Trying to better understand the Riemann hypothesis, I have developed
during the past few years my own ideas via experiments. I am now suspecting
that the problem ($RH$) can neither succumb to an analytic approach
only nor to an arithmetic approach only. A mixture of arithmetic and
analysis seems required and I think to have succeeded in capturing
this duality via good variation theory which depends both on arithmetic
and on complex analysis. 

\begin{flushleft}
Good variation theory emerged in 2010 when I came across the following
two facts
\par\end{flushleft}

\begin{equation}
\sum_{k=1}^{n}\lambda(k)\left\lfloor \frac{n}{k}\right\rfloor =\left\lfloor \sqrt{n}\right\rfloor \end{equation}
which is easy to prove, where $\lambda(k)=(-1)^{\Omega(k)}$ is the
Liouville lambda function, and 

\begin{equation}
\sum_{k=1}^{n}\lambda(k)\ll n^{1/2+\varepsilon}\end{equation}
which is a statement equivalent to the Riemann hypothesis \cite{key-1-1}.
The appearence of the square root of $n$ in the RHS of (1) and (2)
is striking so that I asked myself naively whether we could have $(1)\Rightarrow(2)$?
It turned out that there was no immediate answer to this question
and I tried to figure out what was going on. The first progress I
made arose when I came across the following tauberian theorem of Ingham
(\cite{key-4-2-1}, \cite{key-4-2}).

\[
na(n)\geq-C\ \wedge\ \lim_{n\rightarrow\infty}\sum_{k=1}^{n}a(k)\frac{k}{n}\left\lfloor \frac{n}{k}\right\rfloor =\ell\Rightarrow\sum_{k=1}^{\infty}a(k)=\ell\]

Indeed, some experiments led me to formulate the following general
Tauberian conjecture. Assuming that the condition $a(n)=O\left(n^{-1}\right)$
is satisfied I claim that we have, letting $\Phi(x):=x\left\lfloor \frac{1}{x}\right\rfloor $
denote the Ingham function

\begin{equation}
A_{\Phi}(n)\sim n^{-\beta}\Rightarrow\begin{cases}
0<\beta<\frac{1}{2} & A(n)\sim\frac{1-\beta^{-1}}{\zeta\left(1-\beta\right)}n^{-\beta}\\
\\\beta\geq\frac{1}{2} & A(n)\ll n^{-1/2+\varepsilon}\end{cases}\ \left(n\rightarrow\infty\right)\end{equation}

\begin{flushleft}
which includes $(1)\Rightarrow(2)$ letting $a(n)=\frac{\lambda(n)}{n}$.
Afterwards it was interesting to introduce the broader concept of
functions of good variation as follows. 
\par\end{flushleft}

\subsubsection*{Function of good variation (primary definition)}

A bounded function $g$ which is Riemann integrable on $]0,1]$ is
a function of good variation ($FGV$) of index $\alpha\left(g\right)$
if considering the discrete Volterra equation of linear type %
\footnote{It is an equation of the type $x(n)=f(n)+\sum_{j=0}^{n}y(n,j)x(j)$
$\left(n\geq0\right)$ where $y(.,.)$ and $f(.)$ are known functions
and $x(.)$ is unknown. Cf. for instance \cite{key-101} for results
and other references related to these equations. %
} $A_{g}(n)=n^{-\beta}$ we have the following two Tauberian properties
\begin{enumerate}
\item $\beta<\alpha(g)\Rightarrow A(n)\sim C(\beta)n^{-\beta}\ \left(n\rightarrow\infty\right)$
where $C(\beta)\neq0$
\item $\beta\geq\alpha(g)\Rightarrow A(n)\ll n^{-\alpha(g)}L(n)$ %
\footnote{$L$ denotes a (Karamata) slowly varying function i.e. $\forall x>0\ \lim_{t\rightarrow\infty}\frac{L(tx)}{L(t)}=1$.%
}
\end{enumerate}
So the conjecture (3) is not directly obtained from this primary definition
where we consider $A_{g}(n)=n^{-\beta}$ not the stronger condition
$A_{g}(n)\sim n^{-\beta}$. For our purpose, however, this doesn't
matter because this primary definition will suffice to formulate a
conjecture interesting enough to write down a new equivalence of $RH$. 

In section 1, I prove that $FGV$ exist taking the simplest ones and
I formulate a general conjecture for continuous functions. 

In section 2, the Ingham function $\Phi$ is generalised to a wider
class of similar discontinuous functions: the $BHF$ (broken harmonic
functions). This allows me to formulate a Tauberian conjecture for
$BHF$ in section 3.

Next in section 4, I prove that the Ingham function satisfies an important
condition of the Tauberian conjecture (the Hardy-Littlewood-Ramanujan
criterion) so that $RH$ would be true. In section 5, I prove that
generalised Ingham functions satisfy this important condition as well
so that the Generalized Riemann Hypothesis would be true. In section
6, I extend the method to a set of $L$ functions with multiplicative
coefficients so that the Grand Riemann Hypothesis would be true.

\section{Functions of good variation exist}

It seems important to exhibit examples of $FGV$ since it is not obvious
to see that they exist. Hereafter, I show that $FGV$ is a consistent
concept by proving that affine functions are $FGV$. In particular,
the function $g(x)=\frac{x+1}{2}$ is a $FGV$ of index $\frac{1}{2}$
which is not self evident at first glance.

\subsection{Theorem}

Let $g$ be the affine function $g(x)=c_{1}x+c_{0}$ where $c_{0},c_{1}>0$.
Then $g$ is a $FGV$ of index $\alpha\left(g\right)=\frac{c_{0}}{c_{1}+c_{0}}$
according to the primary definition of $FGV$. 

\begin{flushleft}
More precisely if $A_{g}(n)=n^{-\beta}$ we have 7 cases to consider
summarized in the following table where $g^{\star}(z)=\frac{c_{1}}{1-z}-\frac{c_{0}}{z}$
is the little Mellin transform of $g$. 
\par\end{flushleft}

\begin{onehalfspace}
\begin{center}
\begin{tabular}{|c||>{\centering}m{0.3\paperwidth}|}
\hline 
\noalign{\vskip\doublerulesep}
Condition on $\beta$ & \centering{}$A(n)$ (as $n\rightarrow\infty$)\tabularnewline
\hline
\noalign{\vskip\doublerulesep}
\hline 
\noalign{\vskip\doublerulesep}
$\beta<\alpha(g)-1$ & \centering{}$-\frac{n^{-\beta}}{\beta g^{\star}\left(\beta\right)}+O\left(n^{-1-\beta}\right)$\tabularnewline
\hline
\noalign{\vskip\doublerulesep}
\hline 
\noalign{\vskip\doublerulesep}
$\beta=\alpha(g)-1$ & \centering{}$-\frac{n^{-\beta}}{\beta g^{\star}\left(\beta\right)}+O\left(n^{-\alpha(g)}\log n\right)$\tabularnewline
\hline
\noalign{\vskip\doublerulesep}
\hline 
\noalign{\vskip\doublerulesep}
$\alpha(g)-1<\beta<0$ & \centering{}$-\frac{n^{-\beta}}{\beta g^{\star}\left(\beta\right)}+O\left(n^{-\alpha(g)}\right)$\tabularnewline
\hline
\noalign{\vskip\doublerulesep}
\hline 
\noalign{\vskip\doublerulesep}
$\beta=0$ & \centering{}$\frac{1}{c_{0}}+O\left(n^{-\alpha(g)}\right)$\tabularnewline
\hline
\noalign{\vskip\doublerulesep}
\hline 
\noalign{\vskip\doublerulesep}
$0<\beta<\alpha(g)$ & \centering{}$-\frac{n^{-\beta}}{\beta g^{\star}\left(\beta\right)}+O\left(n^{-\alpha(g)}\right)$\tabularnewline
\hline
\noalign{\vskip\doublerulesep}
\hline 
\noalign{\vskip\doublerulesep}
$\beta=\alpha(g)$ & \centering{}$\left(1-\alpha(g)\right)n^{-\alpha(g)}\log n+O\left(n^{-\alpha(g)}\right)$\tabularnewline
\hline
\noalign{\vskip\doublerulesep}
\hline 
\noalign{\vskip\doublerulesep}
$\beta>\alpha(g)$ & \centering{}$O\left(n^{-\alpha(g)}\right)$\tabularnewline
\hline
\noalign{\vskip\doublerulesep}
\end{tabular}
\par\end{center}
\end{onehalfspace}

\subsection{Proof of theorem 1.1}

Without loss of generality, we take $c_{1}=c_{0}=1/2$ so that $\alpha\left(g\right)=\frac{1}{2}$
and prove the formula of theorem 1.1. for the case $\beta=\alpha(g)-1=-1/2$.
The other formulas are proved similarly for any $c_{1},c_{0}>0$ and
any $\beta$. So let 
\begin{itemize}
\item $g(x)=\frac{x}{2}+\frac{1}{2}$ 
\item $A_{g}(n)=n^{1/2}$ 
\item $h(n)=n^{-1}\left(n^{3/2}-(n-1)^{3/2}\right)$ 
\end{itemize}
First we get the exact formula for any $n\geq2$ (details are omitted)

\begin{equation}
A_{g}(n)=n^{1/2}\Rightarrow A(n)=h(n)+\frac{\left(1/2\right)_{n}}{n!}\left(2+\sum_{k=2}^{n-1}h(k)\frac{k!}{\left(1/2\right)_{k}}\right)\end{equation}

\begin{flushleft}
where $\left(x\right)_{n}=x(x+1)...(x+n-1)$ and it is easy to see
that we have the 3 asymptotic formulas as $n\rightarrow\infty$ 
\par\end{flushleft}

\[
h(n)=\frac{3}{2}n^{-1/2}-\frac{3}{8}n^{-3/2}+O\left(n^{-5/2}\right)\]

\[
\Gamma\left(1/2\right)\frac{(1/2)_{n}}{n!}=n^{-1/2}-\frac{1}{8}n^{-3/2}+O\left(n^{-5/2}\right)\]

\[
\frac{1}{\Gamma\left(1/2\right)}h(k)\frac{k!}{(1/2)_{k}}=\frac{3}{2}-\frac{3}{16}n^{-1}+O\left(n^{-2}\right)\]

\begin{flushleft}
hence plugging these 3 asymptotic formulas in (4) we get 
\par\end{flushleft}

\[
A_{g}(n)=n^{1/2}\Rightarrow A(n)=\frac{3}{2}n^{1/2}-\frac{3}{16}n^{-1/2}\log n+O\left(n^{-1/2}\right)\]

\begin{flushright}
$\square$
\par\end{flushright}

\subsection{A conjecture for continuous function}

Experiments show that much more is true. Indeed, the following conjecture
is very well supported by experiments.

\subsubsection*{Conjecture}

Let $g$ be continuous on $\left[0,1\right]$ satisfying $g(0)g(1)\neq0$.
Then $g$ is a $FGV$ of index $\alpha(g)=\eta(g)$ according to the
primary definition of $FGV$ and supposing that $A_{g}(n)=n^{-\beta}$
we have:
\begin{itemize}
\item $\beta<\alpha(g)\Rightarrow A(n)\sim\left(-\frac{1}{\beta g^{\star}\left(\beta\right)}\right)n^{-\beta}\ \left(n\rightarrow\infty\right)$ 
\item $\beta\geq\alpha(g)\Rightarrow A(n)\ll n^{-\alpha(g)}L(n)$, 
\end{itemize}
where $L$ is slowly varying.

\subsubsection*{Remark}

I didn't work out all the details but I think that the proof could
be derived from the fact that polynomials are $FGV$ and using the
Weierstrass approximation theorem.

\section{Broken harmonic functions}

Almost all tested bounded continuous and discontinuous functions on
$]0,1]$ seem to be $FGV$ and often the conjecture 1.3 works yielding
$\alpha(g)=\eta(g)$ but it is not always the case (see for instance
\cite{key-4} and example 20.1). Anyhow what makes good variation
interesting from a number theoric view point relies on the Ingham
function. Hence the quest for a better and deeper understanding of
the problem led me to consider a set of functions sharing the main
properties of the Ingham function. Some thoughts and many experiments
led me to the natural choice of the so called broken harmonic functions
($BHF$). I will only consider these specific functions in the sequel.

\subsection{Definition of $BHF$}

A bounded function $g$ is a $BHF$ if it exists a real positive sequence
$\left(u_{n}\right)_{n\geq1}$ satisfying $1=u_{1}>u_{2}>u_{3}>...>u_{\infty}=0$
and such that for any $n\geq1$ we have $u_{n+1}<x\leq u_{n}\Rightarrow g(x)=v_{n}x$
where $v_{n}>0$ is an increasing sequence of reals such that $\forall i\geq1$
we have $u_{i}v_{i}<M$ for a constant $M>0$.

\paragraph*{Examples of $BHF$}
\begin{itemize}
\item The Ingham function $\Phi(x)=x\left\lfloor \frac{1}{x}\right\rfloor $
for which $u_{i}=\frac{1}{i}$ and $v_{i}=i$. 
\item For $\lambda>1$, $g_{\lambda}(x)=x\lambda^{\left\lfloor -\frac{\log x}{\log\lambda}\right\rfloor }$
for which $u_{i}=\frac{1}{\lambda^{i-1}}$ and $v_{i}=\lambda^{i-1}$. 
\end{itemize}

\subsection{The $HLR$ criterion}

A function $g$ satisfies the $HLR$ criterion if for any $\beta\geq0$
we have the property \[
A_{g}(n)=n^{-\beta}\Rightarrow\forall\varepsilon>0\ \lim_{n\rightarrow\infty}a(n)n^{1-\varepsilon}=0\]

This criterion was discovered after studying the functions $g_{\lambda}(x)=x\lambda^{\left\lfloor -\frac{\log x}{\log\lambda}\right\rfloor }$
which are $BHF$ satisfying the $HLR$ criterion and the extended
conjecture 1.3 if and only if $\lambda\geq2$ is an integer (for more
details see \cite{key-4}). The name comes from the Hardy-Littlewood
condition in their first Tauberian theorems and from a conjecture
of Ramanujan on the size of coefficients of Dirichlet series in the
Selberg class. This criterion acts as a bridge between arithmetic
and complex analysis and the following Tauberian conjecture illustrates
qualitatively this connection. In this note, I won't formulate the
quantitative Tauberian conjectures described in \cite{key-4} since
the goal is too keep this presentation as short as possible. 

In the sequel {}``$g$ is $HLR$'' means that $g$ satisfies the
$HLR$ criterion and we will often use the following lemma.

\subsection{An important lemma}

If $\left(u_{n}\right)_{n\geq1}$ is a sequence satisfying $u_{1}\neq0$,
then the little Mellin transform of the function $g_{u}$ defined
by 

\[
g_{u}(x)=x\sum_{1\leq k\leq x^{-1}}u_{k}\left\lfloor \frac{1}{kx}\right\rfloor \]
 which is a $BHF$ (details omitted) is given by

\[
g_{u}^{\star}(z)=\frac{\zeta(1-z)U(1-z)}{1-z}\]
 where $U(s)$ is the analytic continuation of the Dirichlet series
$\sum_{n\geq1}u_{n}n^{-s}$.

\subsubsection*{Proof}

By definition for $\Re z<0$ we have

\[
g_{u}^{\star}(z)=\int_{0}^{1}g_{u}(t)t^{-z-1}dt=\sum_{n\geq1}\int_{(n+1)^{-1}}^{n^{-1}}\left(\sum_{1\leq k\leq x^{-1}}u_{k}\left\lfloor \frac{1}{kt}\right\rfloor \right)t^{-z}dt\]

\begin{flushleft}
hence making the variable change $t=x^{-1}$ and letting $w(n)=\sum_{1\leq k\leq n}u_{k}\left\lfloor \frac{n}{k}\right\rfloor $
we get 
\par\end{flushleft}

\[
g_{u}^{\star}(z)=\frac{1}{1-z}\sum_{n\geq1}\frac{w(n)-w(n-1)}{n^{1-z}}\]
 and since we have also $w(n)=\sum_{k=1}^{n}\sum_{d\mid k}u_{d}$
the above equality becomes

\[
g_{u}^{\star}(z)=\frac{1}{1-z}\sum_{n\geq1}\frac{\sum_{d\mid n}u_{d}}{n^{1-z}}=\frac{1}{1-z}\sum_{n\geq1}\frac{(u\star1)(n)}{n^{1-z}}=\frac{\zeta(1-z)U(1-z)}{1-z}\]

\begin{flushright}
$\square$ 
\par\end{flushright}

\section{The anti $HLR$ conjecture for $BHF$ }

Let $g$ be a $BHF$ such that:
\begin{itemize}
\item $\lim_{x\rightarrow0}g(x)\neq0$ exists
\item $(1-z)g^{\star}(z)$ satisfies a Riemann functional equation
\end{itemize}
Then if $g^{\star}$ has a zero in the half-plane $\Re z<\frac{1}{2}$
$g$ is not $HLR$.

\subsection{Corrolary of the conjecture 3}

Suppose $g$ is as above and suppose that $g$ is $HLR$. Then the
non trivial zeros of $g^{\star}$ are on the critical line.

\paragraph*{Proof of corrolary 3.1}

It is simply the contrapositive statement of the anti $HLR$ conjecture.

\subsection{Illustration of the anti $HLR$ conjecture for $BHF$ }

\subsubsection*{The Heilbronn-Davenport zeta function}

The Heilbronn-Davenport zeta function is a good example illustrating
the anti $HLR$ conjecture. Let $\xi=\frac{-2+\sqrt{10-2\sqrt{5}}}{\sqrt{5}-1}=0.284079...$then
Davenport and Heilbronn \cite{key-1} considered the analytic continuation
of the Dirichlet series $H(s)=\sum_{n\geq1}\frac{h(n)}{n^{s}}$ where
$h$ is the 5-periodic sequence $1,\xi,-\xi,-1,0,...$ and proved
that it has nontrivial zeros off the critical line despite the fact
that $H$ satisfies a Riemann functional equation. Then considering
the $BHF$

\[
g_{H}(x):=x\sum_{1\leq k\leq\left\lfloor \frac{1}{x}\right\rfloor }h(k)\left\lfloor \frac{1}{kx}\right\rfloor \]
we have $g_{H}^{\star}(z)=\frac{\zeta(1-z)H(1-z)}{1-z}$ (cf. lemma
2.3) and experiments show clearly that $g_{H}$ isn't $HLR$ (cf.
Fig.1 in the concluding remarks). Therefore the conjecture 3 works
since $g_{H}^{\star}$ has zeros in the half plane $\Re z<\frac{1}{2}$.

\section{The Riemann hypothesis}

Here I prove that the Ingham function is $HLR$ and the reader will
see that it is pure arithmetic involving the fundamental theorem of
arithmetic. Then the Riemann hypothesis follows.

\subsection{Theorem}

The Ingham function $\Phi$ is $HLR$ and more precisely for any $\beta\geq0$
we have 

\[
A_{\Phi}(n)=n^{-\beta}\Rightarrow\ a(n)=O(n^{-1})\]

\subsubsection*{Remark}

Although the Ramanujan condition (the extra $\varepsilon$) is not
necessary for $\Phi$ it will be necessary for the generalised Ingham
functions (see section 5). In some way this means that Hardy-Littlewood
Tauberian condition is sufficient for $RH$ but one needs the deepness
of Ramanujan conjecture to generalise $RH$.

\subsection{Proof of the theorem 4.1 }

In order to prove the theorem I need a lemma.

\subsubsection{Lemma}

Suppose that $f$ is a multiplicative function satisfying $0<f(n)\leq1$
for $n\geq1$. Then we have $\sum_{d\mid n}\mu\left(\frac{n}{d}\right)f(d)=O(1)$.

\subsubsection*{Proof of lemma 4.2.1}

It is easy to see that $b(n):=\sum_{d\mid n}\mu\left(\frac{n}{d}\right)f(d)$
is the multiplicative function given by $b_{p^{v}}=f\left(p^{v}\right)-f\left(p^{v-1}\right)$.
Next letting $n=\prod p_{i}^{\alpha_{i}}$ where $p_{i}$ are the
distinct primes in the factorisation of $n$ we have $-1\leq f\left(p_{i}^{\alpha_{i}}\right)-f\left(p_{i}^{\alpha_{i}-1}\right)\leq1$
hence we get

\[
b_{n}=\prod\left(f\left(p_{i}^{\alpha_{i}}\right)-f\left(p_{i}^{\alpha_{i}-1}\right)\right)=O(1)\]

\begin{flushright}
$\square$ 
\par\end{flushright}

\subsubsection{Proof of the theorem 4.1}

The theorem is true for the case $\beta=0$ since we have trivially
in this case $a_{1}=1$ and $a_{n}=0$ for $n\geq2$. So I consider
$\beta>0$. It is well known that we have

\[
\sum_{k=1}^{n}ka_{k}\left\lfloor \frac{n}{k}\right\rfloor =\sum_{k=1}^{n}\sum_{d\mid k}da_{d}\Rightarrow\ \sum_{d\mid n}da_{d}=n^{1-\beta}-(n-1)^{1-\beta}\]

\begin{flushleft}
Therefore by Möbius inversion we get 
\par\end{flushleft}

\begin{equation}
na_{n}=\sum_{d\mid n}\mu\left(\frac{n}{d}\right)\left(d^{1-\beta}-(d-1)^{1-\beta}\right)\end{equation}

\begin{flushleft}
next we have $d^{1-\beta}-(d-1)^{1-\beta}=(1-\beta)d^{-\beta}+O\left(d^{-1-\beta}\right)$
thus (5) becomes
\par\end{flushleft}

\begin{equation}
na_{n}=(1-\beta)\sum_{d\mid n}\mu\left(\frac{n}{d}\right)d^{-\beta}+\sum_{d\mid n}\mu\left(\frac{n}{d}\right)O(d^{-1-\beta})\end{equation}
Now since $\beta>0$ we have, on one hand, from the lemma 4.2.1 $\sum_{d\mid n}\mu\left(\frac{n}{d}\right)d^{-\beta}=O(1)$
and, on the other hand, $\sum_{n\geq1}n^{-1-\beta}$ converges toward
the finite value $\zeta(1+\beta)$. Consequently we get $\sum_{d\mid n}\mu\left(\frac{n}{d}\right)O(d^{-1-\beta})=O(1)$.

\begin{flushleft}
As a result, (6) becomes $na_{n}=O(1)$ and $\Phi$ is $HLR$. 
\par\end{flushleft}

\begin{flushright}
$\square$ 
\par\end{flushright}

\subsection{Corrolary}

The conjecture 3.1 and the theorem 4.1 imply that $RH$ is true.

\subsubsection*{Proof}

We have
\begin{itemize}
\item $\lim_{x\rightarrow0}\Phi(x)=1\neq0$ 
\item from lemma 2.3,  $(1-z)\Phi^{\star}(z)=\zeta(1-z)$ satisfies a Riemann
functional equation
\item $\Phi$ is $HLR$ from theorem 4.1
\end{itemize}
hence the conjecture 3.1. tells us that $RH$ is true for $\zeta(s)$. 

\begin{flushright}
$\square$
\par\end{flushright}

\section{The Generalized Riemann Hypothesis}

In this section I prove that some generalised Ingham functions are
$HLR$ so that the anti $HLR$ conjecture implies that the Generalized
Riemann Hypothesis is true.

\subsection{Theorem 5.1}

The generalised Ingham function $\Phi_{\chi}(x)=x\sum_{1\leq k\leq1/x}\chi(k)\left\lfloor \frac{1}{kx}\right\rfloor $
is $HLR$ where $\chi$ is a Dirichlet character.

\subsection{Proof of the theorem 5.1 }

The proof is based on the following two lemmas (details are omitted):

\subsubsection{Lemma }

Suppose that $w$ is a completely multiplicative function then $w$
has a Dirichlet inverse given by $w{}^{-1}(n)=\mu(n)w(n)$.

\subsubsection{Lemma }

If $f(n)$ is defined for $n\geq1$ by

\[
f(n)=\sum_{k=1}^{n}u(k)\sum_{1\leq i\leq n/k}v(i)\left\lfloor \frac{n}{ik}\right\rfloor \]

\begin{flushleft}
then we have with the convention $f(0)=0$ 
\par\end{flushleft}

\[
\sum_{d\mid n}u\left(d\right)v\left(\frac{n}{d}\right)=\sum_{d\mid n}\mu\left(\frac{n}{d}\right)\left(f(d)-f(d-1)\right)\]

\subsubsection{Proof of the theorem 5.1}

From lemma 5.2.2 we have

\begin{equation}
A_{\Phi_{\chi}}(n)=n^{-\beta}\Rightarrow\sum_{d\mid n}da\left(d\right)\chi\left(\frac{n}{d}\right)=\sum_{d\mid n}\mu\left(\frac{n}{d}\right)\left(d^{1-\beta}-(d-1)^{1-\beta}\right)\end{equation}

\begin{flushleft}
From 4.2.2. we know that for $\beta\geq0$ we have 
\par\end{flushleft}

\[
\sum_{d\mid n}\mu\left(\frac{n}{d}\right)\left(d^{1-\beta}-(d-1)^{1-\beta}\right)=O(1)\]
Hence letting $a'(n)=na(n)$ and $b(n)=\sum_{d\mid n}\mu\left(\frac{n}{d}\right)\left(d^{1-\beta}-(d-1)^{1-\beta}\right)$
we get from (7) and lemma 5.2.1 (Dirichlet characters are completely
multiplicative)

\[
a'\star\chi(n)=b(n)\Rightarrow a'(n)=b\star\chi^{-1}(n)=\sum_{d\mid n}b\left(d\right)\chi\left(\frac{n}{d}\right)\mu\left(\frac{n}{d}\right)\]

\begin{flushleft}
whence since $\left|b(n)\right|\leq C$ (for some $C>0$) we get 
\par\end{flushleft}

\[
\left|a'(n)\right|\leq\sum_{d\mid n}\left|b\left(d\right)\right|\left|\chi\left(\frac{n}{d}\right)\right|\left|\mu\left(\frac{n}{d}\right)\right|\leq C\tau(n)\ll n^{\varepsilon}\]

\begin{flushleft}
so that $na(n)=O(n^{\varepsilon})$ and $\Phi_{\chi}$ is $HLR$. 
\par\end{flushleft}

\begin{flushright}
$\square$ 
\par\end{flushright}

\subsubsection*{Remark }

Here we can't have something like the Ingham function; i.e., $na(n)=O(1)$.
Indeed let us consider $\chi(n)=1,0,-1,0,1,0,-1,0,...$ and 

\[
\sum_{k=1}^{n}a(k)\Phi_{\chi}\left(\frac{k}{n}\right)=n^{-1}\]

\begin{flushleft}
Let $\left(p_{1}(n)\right)_{n\geq1}$ denotes the increasing sequence
of primes congruent to $1$ modulo $4$; i.e., $p_{1}(1)=5,p_{1}(2)=13,p_{1}(3)=17,p_{1}(4)=29$,
etc. Letting $P(m)=\prod_{i=1}^{m}p_{1}(i)$ it can be shown that
we have for any integer value $m\geq1$ 
\par\end{flushleft}

\[
\left|P(m)a\left(P(m)\right)\right|=2^{m}\]

\begin{flushleft}
therefore $na(n)$ is unbounded. 
\par\end{flushleft}

\subsection{Corrolary}

The conjecture 3.1 and the theorem 5.1 imply that $RH$ is true for
$L(s,\chi)$ where $\chi$ is a Dirichlet character.

\subsubsection*{Proof of corrolary 5.3}

We have
\begin{itemize}
\item $\lim_{x\rightarrow0}\Phi_{\chi}(x)=L\left(1,\chi\right)\neq0$ 
\item from lemma 2.3, $(1-z)\Phi_{\chi}^{\star}(z)=\zeta(1-z)L(1-z,\chi)$
satisfies a Riemann functional equation
\item $\Phi_{\chi}$ is $HLR$ from theorem 5.1
\end{itemize}
hence, from the conjecture 3.1, $RH$ is true for both $\zeta(s)$
and $L(s,\chi)$. 

\begin{flushright}
$\square$ 
\par\end{flushright}

\section{The Grand Riemann Hypothesis}

The previous method extends naturally to the Grand Riemann hypothesis.
In order to do this I need to prove the following theorem.

\subsection{Theorem }

The generalised Ingham function $\Phi_{u}(x)=x\sum_{1\leq k\leq1/x}u(k)\left\lfloor \frac{1}{kx}\right\rfloor $
is $HLR$ where $u$ is any multiplicative function satisfying the
Ramanujan condition $u(n)=O\left(n^{\varepsilon}\right)$.

\subsection{Proof of theorem 6.1}

Before proving this theorem two lemmas are in order. I will prove
only the lemma 6.2.2. The lemma 6.2.1 is classical.

\subsubsection{Lemma}

If $u$ is multiplicative then the Dirichlet inverse $u^{-1}$ is
also multiplicative.

\subsubsection{Lemma}

If $u$ is multiplicative with $u(n)=O(n^{\varepsilon})$ then its
Dirichlet inverse $u^{-1}$ satisfies also $u^{-1}(n)=O\left(n^{\varepsilon}\right)$.

\subsubsection*{Proof of lemma 6.2.2}

Let us fix $u(1)=1$ so that $u^{-1}(1)=1$ then it is well known
that we have the recursive formula for $n\geq2$

\[
u^{-1}(n)=-\sum_{d\mid n,d<n}u^{-1}(d)u\left(\frac{n}{d}\right)\]
Since $u$ is multiplicative so does $u^{-1}$ from lemma 6.2.1. Hence
it suffices to evaluate $u^{-1}\left(p^{n}\right)$ for $p\geq2$
prime and $n\geq1$. From the recursive formula we get 

\begin{equation}
u^{-1}\left(p^{n}\right)=-\sum_{k=0}^{n-1}u^{-1}\left(p^{k}\right)u\left(p^{n-k}\right)\end{equation}

\begin{flushleft}
We assume that $\forall\varepsilon>0,\ u\left(p^{n-k}\right)=O\left(p^{\varepsilon(n-k)}\right)$.
Suppose also that we have the reccurence hypothesis
\par\end{flushleft}
\begin{itemize}
\item $\forall\varepsilon>0,\ u^{-1}\left(p^{k}\right)=O\left(p^{\varepsilon k}\right)$
for any $p\geq2$ prime and any $k\leq n-1$ 
\end{itemize}
Then (8) becomes

\[
\forall\varepsilon>0,\left|u^{-1}\left(p^{n}\right)\right|\leq\sum_{k=0}^{n-1}\left|u^{-1}\left(p^{k}\right)\right|\left|u\left(p^{n-k}\right)\right|\ll\sum_{k=0}^{n-1}p^{\frac{\varepsilon}{2}k}p^{\frac{\varepsilon}{2}(n-k)}=np^{\frac{\varepsilon}{2}n}\]
 next $n\ll p^{\frac{\varepsilon}{2}n}$ for any $\varepsilon>0$
and any $p\geq2$ hence we get 

\[
\left|u^{-1}\left(p^{n}\right)\right|\ll p^{\varepsilon n}\]

\begin{flushleft}
and the reccurence hypothesis is true for all $n$. Thus letting $n=\prod p_{i}^{\alpha_{i}}$
we get 
\par\end{flushleft}

\[
\left|u^{-1}\left(n\right)\right|\ll\prod p_{i}^{\varepsilon\alpha_{i}}=n^{\varepsilon}\]

\begin{flushright}
$\square$ 
\par\end{flushright}

\subsubsection{Proof of theorem 6.1}

\begin{flushleft}
From lemma 5.2.2 we have
\par\end{flushleft}

\begin{equation}
A_{\Phi_{u}}(n)=n^{-\beta}\Rightarrow\sum_{d\mid n}da\left(d\right)u\left(\frac{n}{d}\right)=\sum_{d\mid n}\mu\left(\frac{n}{d}\right)\left(d^{1-\beta}-(d-1)^{1-\beta}\right)\end{equation}

\begin{flushleft}
From 4.2.2. we know that for $\beta\geq0$ we have 
\par\end{flushleft}

\[
\sum_{d\mid n}\mu\left(\frac{n}{d}\right)\left(d^{1-\beta}-(d-1)^{1-\beta}\right)=O(1)\]

\begin{flushleft}
Hence letting $a'(n)=na(n)$ and $b(n)=\sum_{d\mid n}\mu\left(\frac{n}{d}\right)\left(d^{1-\beta}-(d-1)^{1-\beta}\right)$
we get from (9) 
\par\end{flushleft}

\[
a'\star u(n)=b(n)\Rightarrow a'(n)=b\star u^{-1}(n)=\sum_{d\mid n}b\left(n/d\right)u^{-1}(d)\]

\begin{flushleft}
whence since $b(n)=O(1)$ and $u^{-1}(n)=O\left(n^{\varepsilon}\right)$
from lemma 6.2.2 we get for any $\varepsilon>0$
\par\end{flushleft}

\[
\left|na(n)\right|\leq\sum_{d\mid n}\left|b\left(n/d\right)\right|\left|u^{-1}\left(d\right)\right|\ll\sum_{d\mid n}d^{\varepsilon/2}\ll n^{\varepsilon/2}\tau(n)\ll n^{\varepsilon}\]

\begin{flushleft}
since $\tau(n)\ll n^{\varepsilon/2}$, consequently $na(n)=O(n^{\varepsilon})$
and $\Phi_{u}$ is $HLR$. 
\par\end{flushleft}

\begin{flushright}
$\square$
\par\end{flushright}

\subsection{Corrolary}

Suppose that:
\begin{itemize}
\item $u$ is multiplicative with $u(n)=O(n^{\varepsilon})$ 
\item the analytic continuation of $U(s)=\sum_{n\geq1}\frac{u(n)}{n^{s}}$
satisfies a Riemann functional equation 
\item $ $$\sum_{n\geq1}\frac{u(n)}{n}$ converges toward a non zero limit
\end{itemize}
Then the conjecture 3.1 and the theorem 6.1 imply that $RH$ is true
for $U(s)$.

\subsubsection*{Proof of corrolary 6.3 }

We have
\begin{itemize}
\item $\lim_{x\rightarrow0}\Phi_{u}(x)=\sum_{n\geq1}\frac{u(n)}{n}\neq0$ 
\item from lemma 2.3, $(1-z)\Phi_{u}^{\star}(z)=\zeta(1-z)U(1-z)$ satisfies
a Riemann functional equation
\item $\Phi_{u}(x)=x\sum_{1\leq k\leq1/x}u(k)\left\lfloor \frac{1}{kx}\right\rfloor $
is $HLR$ from theorem 6.1
\end{itemize}
hence, from the conjecture 3.1, $RH$ is true for both $\zeta(s)$
and $U(s)$.

\subsubsection*{Example}

Let $\tau_{r}(n)$ denote the Ramanujan tau numbers (\cite{key-102},
pp. 114, 131) defined by

\[
\sum_{n\geq1}\tau_{r}(n)x^{n}=x\prod_{n\geq1}\left(1-x^{n}\right)^{24}\]
then $u(n)=\frac{\tau_{r}(n)}{n^{11/2}}$ satisfies the conditions
of corrolary 6.3 since letting

\[
U(s)=\sum_{n\geq1}\frac{\tau_{r}(n)}{n^{s+11/2}}\]

\begin{flushleft}
which satisfies a Riemann functional equation we have:
\par\end{flushleft}
\begin{itemize}
\item \begin{flushleft}
$u$ is multiplicative since $\tau_{r}$ is multiplicative (conjectured
by Ramanujan and proved soon after by Mordell in 1917 \cite{key-201})
\par\end{flushleft}
\item \begin{flushleft}
$u(n)=O\left(n^{\varepsilon}\right)$ since $\tau_{r}(n)=O\left(n^{11/2+\varepsilon}\right)$
(conjectured by Ramanujan and proved much later by Deligne in 1974
\cite{key-3})
\par\end{flushleft}
\item \begin{flushleft}
$\lim_{x\rightarrow0}\Phi_{u}(x)=\sum_{n\geq1}\frac{\tau_{r}(n)}{n^{13/2}}=0.8...\neq0$ 
\par\end{flushleft}
\item \begin{flushleft}
$(1-z)\Phi_{u}^{\star}(z)=\zeta(1-z)U(1-z)$ satisfies a Riemann functional
equation
\par\end{flushleft}
\end{itemize}
Hence, from corrolary 6.3, $U(s)$ satisfies the Riemann hypothesis.

\section*{Concluding remarks}

The Davenport-Heilbronn example described in 3.2. is interesting on
its own. It is known (see for instance \cite{key-11}) that the non
trivial zeros of $H(s)$ have real part dense in $]0,1[$ and it could
be a general property. 

Namely I claim that if $F(s)=\sum_{n\geq1}\frac{f(n)}{n^{s}}$ has
an analytic continuation, satisfies a Riemann functional equation
and has some zeros off the critical line then it has in fact infinitely
many zeros off the critical line and the real part of these zeros
are dense in $]0,1[$.

In this case the generalised Ingham function 

\[
g_{f}(x)=x\sum_{k\leq x^{-1}}f(k)\left\lfloor \frac{1}{kx}\right\rfloor \]
isn't $HLR$ from the anti $HLR$ conjecture and I claim that we have
for any $\beta\geq0$ 

\[
A_{g_{f}}(n)=n^{-\beta}\Rightarrow\forall\varepsilon>0\ \lim_{n\rightarrow\infty}a_{n}n^{1/2-\varepsilon}=0\ \wedge\limsup_{n\rightarrow\infty}\left|a_{n}n^{1/2+\varepsilon}\right|=+\infty\]
which seems supported by experiments using the Heilbronn-Davenport
counter-example and the associated $BHF$ $g_{H}$ defined in 3.2.
as shown by the graphic below (Fig. 1) which looks bounded by a slowly
varying function like $\log$.

\begin{center}
Fig.1) Plot of $n^{1/2}a(n)$ where $A_{g_{H}}(n)=n^{-1/3}$ 
\par\end{center}

\begin{center}
\includegraphics[width=0.6\paperwidth]{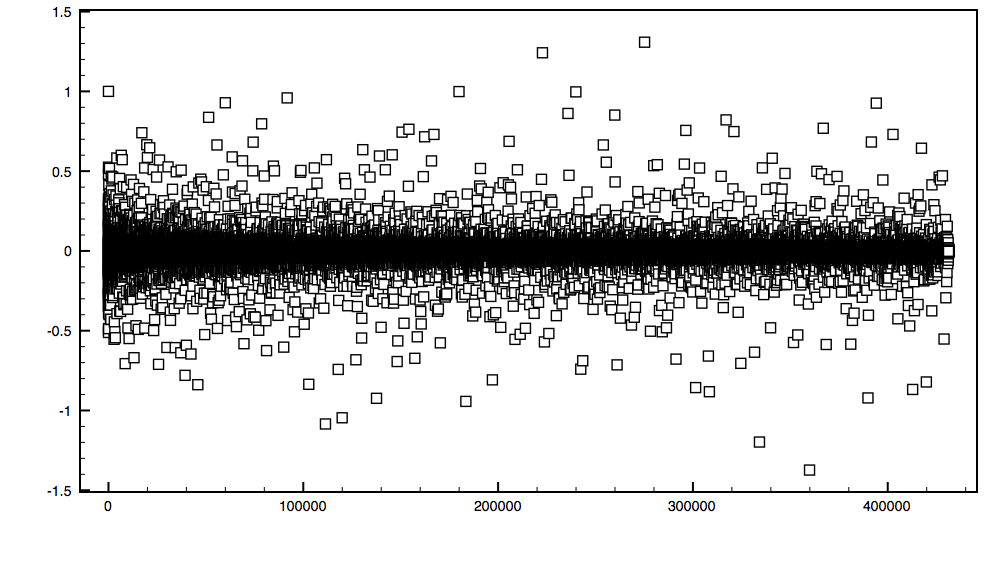}
\par\end{center}

In particular if $\zeta$ has a zero off the critical line it should
have in fact infinitely many zeros off the critical line and we could
find zeta zeros as close as we wish from the line $x=1$. In some
way this is supported by the best zero free regions known to this
day which are always asymptotically close to the line $x=1$ (although
it is not the latest see \cite{key-105} for recent work on this topic).

Hence one should explore this kind of property of values distribution
of little Mellin transform of $BHF$ like the above function $g_{f}$
in order to better understand properties of nontrivial zeros (simplicity
of zeros, Montgomery's pair correlation conjecture, independence of
imaginary parts over the rationals,...).

To prove the anti $HLR$ conjecture however, I believe that an algebraic
approach could be promising. The idea is to consider the space of
$FGV$ and the subspace of affine functions by parts and to look for
Tauberian invariants. BTW a direct proof of the anti $HLR$ conjecture
looks possible exploring sums over zeros of the little Mellin transform.
For instance if $A_{\Phi}(n)=n^{-\beta}$ one should have this kind
of explicit formula

\[
A(n)=-\frac{n^{-\beta}}{\beta\Phi^{\star}\left(\beta\right)}+\sum_{\rho}c_{\beta}(\rho)n^{-\rho}+\sum_{k\geq1}d_{\beta}(k)n^{-2k}\]
 where $\rho$ denote the nontrivial zeros of $\zeta$ and $c_{\beta}(\rho),$$d_{\beta}(k)$
are suitable coefficients. 

With this approach we link zeros of $L$ functions to multiplicative
number theory. Thus the reader should note the consistency with what
is suggested by experts in analytic number theory, i.e., $RH$ is
true for analytic continuation of Dirichlet series satisfying a Riemann
functional equation if and only if there is an Euler product.

\section*{Acknowledgments}

I warmly thank Doron Zeilberger for his ongoing support during this
research and I thank the referee for his time.

\end{document}